# ESSCHER TRANSFORM AND THE DUALITY PRINCIPLE FOR MULTIDIMENSIONAL SEMIMARTINGALES

By Ernst Eberlein, Antonis Papapantoleon[1]
and Albert N. Shiryaev[2]

*University of Freiburg, Vienna University of Technology and Steklov Mathematical Institute*

The duality principle in option pricing aims at simplifying valuation problems that depend on several variables by associating them to the corresponding dual option pricing problem. Here, we analyze the duality principle for options that depend on several assets. The asset price processes are driven by general semimartingales, and the dual measures are constructed via an Esscher transformation. As an application, we can relate swap and quanto options to standard call and put options. Explicit calculations for jump models are also provided.

**1. Introduction.** In this work, we continue our study of the duality principle in option pricing for general semimartingales, initiated in Eberlein, Papapantoleon and Shiryaev (2008) (henceforth, EPS). Here, we consider options that depend on several assets; the valuation of such options requires the knowledge of the joint distribution or characteristic function, and a high-dimensional integration. We aim to simplify this valuation problem by relating it to its dual option pricing problem.

In one-dimensional semimartingale models, we operate in EPS with a single dual measure and the density process is the asset price itself. In the multidimensional case though, we cannot work with a single measure; indeed, the density of the corresponding dual measure will be dictated by the problem at hand. It turns out that the Esscher change of measure [cf. Kallsen and Shiryaev (2002)] is the appropriate concept to describe the density between the original and the dual measure. Therefore, our main re-

Received September 2008; revised December 2008.

[1]Supported by the Austrian Science Fund FWF Grant Y328, START Prize.

[2]Supported by Grants RFBI 08-01-00740 and 08-01-91205-YaF, and also by a grant of the University of Freiburg.

sult describes the triplet of predictable characteristics of one-dimensional semimartingales—defined as the inner product of a vector with the driving multidimensional semimartingale—under the dual measure.

As an application of our results, we relate swap options and quanto options to standard European call or put options, for general semimartingale models. For semimartingales with independent increments, we can derive a duality relationship between an option depending on three assets and a standard call or put option. This yields a significant reduction in the computational complexity of these valuation problems.

The model we employ to describe the evolution of asset price processes is an exponential semimartingale model. Semimartingales are the most natural and general processes we can consider from the point of view of arbitrage theory; they also contain, as subclasses, most of the models used in mathematical finance, such as Brownian motion and general diffusions, Lévy processes, affine processes, time-changed Lévy models and stochastic volatility models. We identify the driving semimartingale process by its triplet of predictable characteristics [cf. Jacod (1979)].

Duality results in the univariate case have been extensively studied in the literature for several subclasses of semimartingales, see the introduction of EPS for a (nonexhaustive) list of references. However, there seems to be considerably less work in the multivariate case; we mention here the articles of Margrabe (1978), Geman, El Karoui and Rochet (1995) and Gerber and Shiu (1996) for the Black–Scholes model, while Eberlein and Papapantoleon (2005) work with time-inhomogeneous Lévy processes; Fajardo and Mordecki (2006) consider Lévy processes, and also handle American style options. Molchanov and Schmutz (2008) derive analogous results using techniques from convex geometry. Schroder (1999) studied the problem in a semimartingale framework, but did not derive a general representation, for example, the characteristics, under the dual measure.

The paper is organized as follows. In Section 2, we review some results on multidimensional semimartingales, their triplet of predictable characteristics, the Laplace cumulant process and linear transformations of multidimensional semimartingales. In Section 3, we describe the asset price model; Section 4 contains the main result of this work, describing one-dimensional semimartingales under the dual measures. Finally, in Section 5, we describe several applications of the duality principle and present some explicit examples, especially for jump processes.

**2. Semimartingales and their characteristics.** 1. Let $\mathbb{R}^d$ denote the $d$-dimensional Euclidean space. The Euclidean scalar product between two vectors $u, v \in \mathbb{R}^d$ is denoted by $\langle u, v \rangle$ or $u^\top v$, where $u^\top$ denotes the transpose of the vector (or matrix) $u$. The Euclidean norm is denoted by $|\cdot|$, and $\mathbf{e}_i$ denotes the unit vector, where the $i$th entry is 1 and all others zero, that is, $\mathbf{e}_i = (0, \ldots, 1, \ldots, 0)^\top$.

The inner product is extended from real to complex numbers as follows: for $u = (u_k)_{1 \leq k \leq d}$ and $v = (v_k)_{1 \leq k \leq d}$ in $\mathbb{C}^d$, set $\langle u, v \rangle := \sum_{k=1}^d u_k v_k$; therefore, we do not use the Hermitian inner product $\sum_{k=1}^d u_k \overline{v}_k$. Moreover, we denote by $iv := (iv_k)_{1 \leq k \leq d}$.

Let $M_d(\mathbb{R})$ denote the space of real $d \times d$ matrices, and let $\| \cdot \|$ denote the norm on the space of $d \times d$ matrices induced by the Euclidean norm on $\mathbb{R}^d$. In addition, let $M_{nd}(\mathbb{R})$ denote the space of real $n \times d$ matrices, and similarly $\| \cdot \|$ denotes the induced norm on this space. Note that we could equally well work with any vector norm on $\mathbb{R}^d$ and the norms induced by, or consistent with, it on $M_d(\mathbb{R})$ and $M_{nd}(\mathbb{R})$.

Define the set $D := \{x \in \mathbb{R}^d : |x| > 1\}$, hence, $D^c$ denotes the closed unit ball in $\mathbb{R}^d$. The function $h = h(x)$ denotes a *truncation function*, that is, a bounded function with compact support that behaves as $h(x) = x$ around the origin; the *canonical choice* is $h(x) = x 1_{\{|x| \leq 1\}} = x 1_{D^c}(x)$. We assume that $h$ satisfies the *antisymmetry* property $h(-x) = -h(x)$.

REMARK 2.1. The truncation function on $\mathbb{R}^n$, $n \neq d$, will also be denoted by $h(x)$, for $x \in \mathbb{R}^n$; that is, the *argument* will determine the *dimension*.

In general, we follow the notation of Jacod and Shiryaev (2003) (henceforth, JS); any unexplained notation is typically used as in JS.

2. Consider a stochastic basis $\mathscr{B} = (\Omega, \mathcal{F}, \mathbf{F}, P)$ in the sense of JS I.1.2, where $\mathcal{F} = \mathcal{F}_T$, $\mathbf{F} = (\mathcal{F}_t)_{0 \leq t \leq T}$ and $T$ is a finite time horizon. Let $H = (H_t)_{0 \leq t \leq T}$ be a $d$-dimensional general semimartingale, that is, $H = (H^1, \ldots, H^d)^\top$ with $H_0 = 0$. Every semimartingale has the *canonical representation* (cf. JS II.2.34)

$$(2.1) \qquad H = H_0 + B + H^c + h(x) * (\mu - \nu) + (x - h(x)) * \mu$$

or, in detail

$$(2.2) \quad H_t = H_0 + B_t + H_t^c + \int_0^t \int_{\mathbb{R}^d} h(x)\, d(\mu - \nu) + \int_0^t \int_{\mathbb{R}^d} (x - h(x))\, d\mu,$$

where:

(a) $B = (B_t)_{0 \leq t \leq T}$ is an $\mathbb{R}^d$-valued predictable process of bounded variation;
(b) $H^c = (H_t^c)_{0 \leq t \leq T}$ is the continuous martingale part of $H$; $H^c$ has the predictable quadratic characteristic $\langle H^c \rangle = C$, which is a predictable $\mathbb{R}^{d \times d}$-valued process of bounded variation, whose values are nonnegative symmetric matrices;
(c) $\nu = \nu(\omega; dt, dx)$ is a predictable random measure on $[0, T] \times \mathbb{R}^d$; it is the compensator of the random measure of jumps $\mu = \mu(\omega; dt, dx)$ of $H$.

Here, $W * \mu$ denotes the integral process, and $W * (\mu - \nu)$ denotes the stochastic integral with respect to the compensated random measure $\mu - \nu$; cf. JS, Chapter II.

The processes $B$, $C$ and the measure $\nu$ are called the *triplet of predictable characteristics* of the semimartingale $H$ with respect to the probability measure $P$, and will be denoted by

$$\mathbb{T}(H|P) = (B, C, \nu).$$

The characteristics are uniquely defined, up to indistinguishability of course.

In addition, there exists an increasing predictable process $A$, predictable processes $b$, $c$ and a transition kernel $F$ from $(\Omega \times [0,T], \mathcal{P})$ into $(\mathbb{R}^d, \mathcal{B}(\mathbb{R}^d))$ such that

(2.3) $$B = b \cdot A, \qquad C = c \cdot A, \qquad \nu = F \otimes A,$$

or, in detail

(2.4) $$B_t = \int_0^t b_s \, dA_s, \qquad C_t = \int_0^t c_s \, dA_s,$$
$$\nu([0,t] \times E) = \int_0^t \int_E F_s(dx) \, dA_s,$$

where $E \in \mathcal{B}(\mathbb{R}^d)$; cf. JS, Proposition II.2.9.

Every semimartingale $H$ with triplet $\mathbb{T}(H|P) = (B, C, \nu)$ can be associated to a *Laplace cumulant process* denoted by $K = (K_t)_{0 \le t \le T}$, defined via

(2.5) $$K(u) = \langle u, B \rangle + \tfrac{1}{2}\langle u, Cu \rangle + (e^{\langle u, x \rangle} - 1 - \langle u, h(x) \rangle) * \nu.$$

Moreover, we have that $K(u) = \kappa(u) \cdot A$, where

(2.6) $$\kappa(u) = \langle u, b \rangle + \frac{1}{2}\langle u, cu \rangle + \int_{\mathbb{R}^d} (e^{\langle u, x \rangle} - 1 - \langle u, h(x) \rangle) F(dx).$$

The Laplace cumulant process satisfies the following property (cf. Corollary II.2.48 in JS):

(2.7) $$\frac{e^{i\langle u, H \rangle}}{\mathcal{E}(K(iu))} \in \mathcal{M}_{\mathrm{loc}}(P)$$

for all $u \in \mathbb{R}^d$, assuming that $\mathcal{E}(K(iu))$ never vanishes; see Remark 2.3 for sufficient conditions. Here, $\mathcal{E}(\cdot)$ denotes the *stochastic exponential*; cf., for example, JS I.4.61. Formula (2.7) is also called the *martingale version of the Lévy–Khintchine formula for semimartingales*.

Note that given a cumulant process (2.5), satisfying (2.7), we can immediately conclude that the triplet of characteristics for the semimartingale $H$ is given by $(B, C, \nu)$; cf. Corollary II.2.48 in JS.

REMARK 2.2. If the characteristics $(B, C, \nu)$ are *absolutely continuous*, then we can choose the process $A_t = t$. Then we call the triplet $(b, c, F)$ the *differential characteristics* of $H$.

REMARK 2.3. The following diagram of statements holds true:

(2.8) $\qquad\qquad (1) \quad \Rightarrow \quad (2) \quad \Leftrightarrow \quad (3) \quad \Rightarrow \quad (4),$

where:

1. $\mathbb{T}(H|P) = (B, C, \nu)$ is absolutely continuous;
2. $H$ has no fixed times of discontinuity;
3. $H$ is a quasi-left-continuous process;
4. $K$ is a continuous process.

The first statement follows by direct calculations; for the others see I.2.25, II.1.19 and III.7.4 in JS.

3. We consider linear transformations and projections of general semimartingales. The following result, which seems to be in the literature already, describes the triplet of predictable characteristics under such a transformation. We provide a short proof and also study some properties of the resulting process. Analogous results for Lévy and time-inhomogeneous Lévy processes can be found in Sato (1999), Proposition 11.10, and Papapantoleon (2007), Proposition 2.10.

Let $H$ be an $\mathbb{R}^d$-valued semimartingale and let $U$ be an $n \times d$-dimensional real valued matrix; then $UH = U \times H$ is an $\mathbb{R}^n$-valued semimartingale and the following proposition determines the triplet of $UH$.

PROPOSITION 2.4. *Consider an $\mathbb{R}^d$-valued semimartingale $H = (H_t)_{0 \leq t \leq T}$ with triplet $\mathbb{T}(H|P) = (B, C, \nu)$ and a real $n \times d$ matrix $U$ ($U \in M_{nd}(\mathbb{R})$). Then $UH = (UH_t)_{0 \leq t \leq T}$ is an $\mathbb{R}^n$-valued semimartingale with triplet of predictable characteristics $\mathbb{T}(UH|P) = (B^U, C^U, \nu^U)$ of the form*

(2.9) $\qquad B^U = b^U \cdot A, \qquad C^U = c^U \cdot A, \qquad \nu^U = F^U \otimes A,$

*where*

(2.10) $\qquad \begin{aligned} b_s^U &= Ub_s + \int_{\mathbb{R}^d} (h(Ux) - Uh(x)) F_s(dx), \\ c_s^U &= Uc_s U^\top, \\ F_s^U(E) &= \int_{\mathbb{R}^d} 1_E(Ux) F_s(dx), \qquad E \in \mathcal{B}(\mathbb{R}^n \setminus \{0\}). \end{aligned}$

*(Recall that $h$ denotes a generic truncation function; cf. Remark 2.1.)*

PROOF. One could prove the statement by directly calculating the characteristics, see analogous results for the stochastic integral process in Kallsen and Shiryaev (2001), Lemma 3, and EPS (Lemma 3.3). However, using the martingale version of the Lévy–Khintchine formula for semimartingales, a very simple proof can be given.

Indeed (cf. Corollary II.2.48 in JS) to prove the assertion it suffices to show that

$$\frac{e^{i\langle z, UH \rangle}}{\mathcal{E}(K^U(iz))} \in \mathcal{M}_{\mathrm{loc}}(P) \tag{2.11}$$

for any $z \in \mathbb{R}^n$; here $K^U$ denotes the cumulant process associated with the triplet $(B^U, C^U, \nu^U)$.

Since $\mathbb{T}(H|P) = (B, C, \nu)$, we immediately have that for any $z \in \mathbb{R}^n$ holds

$$\frac{e^{i\langle U^\top z, H \rangle}}{\mathcal{E}(K(iU^\top z))} \in \mathcal{M}_{\mathrm{loc}}(P),$$

where

$$\frac{e^{i\langle U^\top z, H \rangle}}{\mathcal{E}(K(iU^\top z))}$$
$$= e^{i\langle U^\top z, H \rangle} / (\mathcal{E}(\langle iU^\top z, B \rangle + 1/2 \langle iU^\top z, CiU^\top z \rangle$$
$$\qquad + (e^{\langle iU^\top z, x \rangle} - 1 - \langle iU^\top z, h(x) \rangle) * \nu))$$
$$= \frac{e^{i\langle z, UH \rangle}}{\mathcal{E}(\langle iz, UB \rangle + 1/2 \langle iz, UCU^\top iz \rangle + (e^{\langle iz, Ux \rangle} - 1 - \langle iz, Uh(x) \rangle) * \nu)}$$
$$= \frac{e^{i\langle z, UH \rangle}}{\mathcal{E}(\langle iz, B^U \rangle + 1/2 \langle iz, UCU^\top iz \rangle + (e^{\langle iz, Ux \rangle} - 1 - \langle iz, h(Ux) \rangle) * \nu)}$$
$$= \frac{e^{i\langle z, UH \rangle}}{\mathcal{E}(K^U(iz))},$$

where $B^U = UB + [h(Ux) - Uh(x)] * \nu$; hence the assertion is proved. □

Next, we derive some results about the properties of the process $UH$.

LEMMA 2.5. *If $H$ is a special semimartingale, then $UH$ is also a special semimartingale.*

PROOF. If suffices to prove that $1_{\{|y|>1\}}|y| * \nu^U \in \mathcal{V}$; cf. JS II.2.29. We have

$$1_{\{|y|>1\}}|y| * \nu^U = 1_{\{|Ux|>1\}}|Ux| * \nu$$

$$\leq 1_{\{|Ux|>1\}}\|U\||x| * \nu$$

(2.12) $$\leq 1_{\{\|U\||x|>1\}}\|U\||x| * \nu$$

(2.13) $$\leq 1_{\{|x|>1\}}\|U\||x| * \nu + 1_{\{1/\|U\|<|x|\leq 1\}}\|U\|^2|x|^2 * \nu \in \mathcal{V},$$

which follows from the assumption that $H$ is a special semimartingale and the properties of the compensator (cf. JS II.2.13). Notice that we have implicitly assumed that $\|U\| \geq 1$; otherwise, we can conclude already from (2.12). □

LEMMA 2.6. *If $H$ is a quasi-left-continuous process, then $UH$ is also a quasi-left-continuous process.*

PROOF. Let $H$ be quasi-left-continuous, then $\nu(\{t\} \times \mathbb{R}^d) = 0$ for all $t \in [0,T]$. Therefore, for the process $UH$, we have that

$$\nu^U(\{t\} \times \mathbb{R}^n) = \int 1_{\mathbb{R}^n \setminus \{0\}}(Ux)\nu(\{t\} \times dx) = 0 \qquad \forall t \in [0,T].$$

Hence, $UH$ is also quasi-left-continuous. □

**3. Exponential semimartingale models.** We present some details about the model we employ, where asset prices are modeled as exponentials of general semimartingales. For the sake of completeness, we also derive the martingale condition in this framework, subject to a mild and natural assumption on the driving processes.

Let $\mathcal{M}_{\text{loc}}(P)$ be the class of all local martingales on the given stochastic basis $(\Omega, \mathcal{F}, (\mathcal{F}_t)_{0\leq t\leq T}, P)$; let $\mathcal{V}$ denote the class of processes with bounded variation. Let $H = (H^1, \ldots, H^d)^\top$ be the vector of semimartingale driving processes; it has the triplet of characteristics $\mathbb{T}(H|P) = (B, C, \nu)$.

ASSUMPTION (𝔼𝕊). Assume that the process $1_{\{|x|>1\}}e^{x^i} * \nu \in \mathcal{V}$ for all $i \in \{1, \ldots, d\}$.

REMARK 3.1. Equivalently, we could assume that the process $H^i = \mathbf{e}_i^\top H$ is *exponentially special*, or that the process $e^{H^i}$ is a *special semimartingale*, for all $i \in \{1, \ldots, d\}$. Additionally, under (𝔼𝕊) the martingale version of the Lévy–Khintchine formula (2.7) holds true for *real* arguments, and in particular for the unit vectors $\mathbf{e}_i$, $i \in \{1, \ldots, d\}$; cf. Proposition II.8.26 and Remark III.7.15 in JS.

The following result further characterizes the set of exponentially special semimartingales; it also extends Theorem 25.17(i) in Sato (1999) to general semimartingales.

LEMMA 3.2. *Let $H$ be an $\mathbb{R}^d$-valued semimartingale with triplet $(B,C,\nu)$. The set $\mathcal{U}$, where*

$$\mathcal{U} = \{u \in \mathbb{R}^d : e^{\langle u,x \rangle} 1_{\{|x|>1\}} * \nu \in \mathcal{V}\},$$

*is a convex set and contains the origin.*

PROOF. The definition of the compensator $\nu$ immediately shows that $\mathcal{U}$ contains the origin (cf. JS II.2.13, I.3.9). Now, consider $u, v \in \mathcal{U}$ and $p, q \in (0,1)$ with $q = 1-p$; using that $F$ in (2.4) is a kernel and applying Hölder's inequality twice, we get (recall that $D = \{x \in \mathbb{R}^d : |x| > 1\}$)

$$
\begin{aligned}
e^{\langle pu+qv,x\rangle} 1_{\{|x|>1\}} * \nu &= \int_0^\cdot \int_D e^{\langle pu+qv,x\rangle} F_s(dx)\, dA_s \\
&\leq \int_0^\cdot \left(\int_D e^{\langle u,x\rangle} F_s(dx)\right)^p \left(\int_D e^{\langle v,x\rangle} F_s(dx)\right)^q dA_s \\
&\leq \left(\int_0^\cdot \int_D e^{\langle u,x\rangle} F_s(dx)\, dA_s\right)^p \left(\int_0^\cdot \int_D e^{\langle v,x\rangle} F_s(dx)\, dA_s\right)^q \\
&= (e^{\langle u,x\rangle} 1_D * \nu)^p (e^{\langle v,x\rangle} 1_D * \nu)^q \in \mathcal{V}.
\end{aligned}
$$

Hence, the set $\mathcal{U}$ is convex.  □

Let $S = (S^1, \ldots, S^d)^\top$ denote the vector of asset price processes. Each component $S^i$ of $S$ is an exponential semimartingale, that is, a stochastic process with representation

$$(3.1) \qquad S_t^i = e^{H_t^i}, \qquad 0 \leq t \leq T, 1 \leq i \leq d,$$

(shortly, $S^i = e^{H^i}$), where $H^i = (H_t^i)_{0 \leq t \leq T}$ is a real-valued semimartingale with canonical representation

$$(3.2) \qquad H^i = H_0^i + B^i + H^{i,c} + h^i(x) * (\mu - \nu) + (x^i - h^i(x)) * \mu,$$

where $h^i(x) = \mathbf{e}_i^\top h(x)$. For simplicity, we assume that $S_0^i = 1$ for all $i \in \{1, \ldots, d\}$; we also assume that the interest rate and dividend yields are zero.

REMARK 3.3. Note that $h^i(x)$ can equally well serve as a truncation function on the real line, instead of $h(x^i)$, $x^i \in \mathbb{R}$. This is a feature specific to the unit vector; for arbitrary vectors $u^\top h(x)$ will look quite different from $h(u^\top x)$.

PROPOSITION 3.4. *Subject to Assumption ($\mathbb{ES}$), we have that*

$$(3.3) \quad S^i = e^{H^i} \in \mathcal{M}_{\mathrm{loc}}(P) \quad \Leftrightarrow \quad B^i + \tfrac{1}{2} C^{ii} + (e^{x^i} - 1 - h^i(x)) * \nu = 0.$$

PROOF. We give two proofs, since they reveal an interesting interplay regarding canonical representations and truncation functions.

A. Consider the unit vector $\mathbf{e}_i^\top$, and apply Proposition 2.4 to this vector and the semimartingale $H$. Then we get that $H^i = \mathbf{e}_i^\top H$ is a real-valued semimartingale with triplet $\mathbb{T}(H^i|P) = (\bar{B}^i, C^i, \nu^i)$, where

$$\bar{b}^i = b^i + [h(x^i) - h^i(x)] * F,$$
$$(3.4) \qquad c^i = c^{ii},$$
$$F^i(E) = F(\{x \in \mathbb{R}^d : x^i \in E\}), \qquad E \in \mathcal{B}(\mathbb{R} \setminus \{0\}).$$

Further, $H^i = (H^i_t)_{0 \leq t \leq T}$ admits the canonical representation

$$(3.5) \qquad H^i = H^i_0 + \bar{B}^i + H^{i,c} + h(y) * (\mu^i - \nu^i) + (y - h(y)) * \mu^i;$$

compare with the representation (3.2).

Now, applying equivalence (3.5) in EPS to the real-valued process $H^i$, we get

$$(3.6) \quad S^i = e^{H^i} \in \mathcal{M}_{\mathrm{loc}}(P) \quad \Leftrightarrow \quad \bar{B}^i + \tfrac{1}{2}C^{ii} + (e^y - 1 - h(y)) * \nu^i = 0,$$

which after some calculations using (3.4) yields the asserted result.

B. We can rewrite equivalence (3.3), using the form of the cumulant process (2.5), as follows:

$$(3.7) \qquad S^i = e^{H^i} \in \mathcal{M}_{\mathrm{loc}}(P) \quad \Leftrightarrow \quad K(\mathbf{e}_i) = 0.$$

Now, the "if" part is rather obvious, using the martingale version of the Lévy–Khintchine formula (2.7) for the real argument $\mathbf{e}_i$.

Conversely, if $e^{H^i} \in \mathcal{M}_{\mathrm{loc}}(P)$, from the uniqueness of the multiplicative decomposition of a special semimartingale [cf. Jacod (1979), VI.2a and Theorem 6.19], we get that

$$(3.8) \qquad \mathcal{E}(K(\mathbf{e}_i)) = 1.$$

Now, we can apply the stochastic logarithm on both sides of (3.8) since $\Delta K(\mathbf{e}_i) > -1$ [cf. Kallsen and Shiryaev (2002), page 405] which leads to the required result. $\square$

**4. Multidimensional dual measures.** The aim of this section is to characterize one-dimensional semimartingales, defined as scalar products of the driving semimartingale $H$ and $d$-dimensional vectors $u$, under a suitable equivalent probability measure. This measure, termed the *dual measure*, is defined by an *Esscher transformation*; cf. Kallsen and Shiryaev (2002) (henceforth, KS; note that we do not use the same notation as KS; in particular, $\widetilde{K}$ in KS is denoted $K$ here and vice versa). We point out that the stochastic integral in the Esscher transform pertains to interest rate modeling.

Let $L(H)$ denote the set of all (predictable) integrable processes $\vartheta$ with respect to the semimartingale $H$ (JS III.6.17); $\vartheta^\top \cdot H$ denotes the stochastic integral of $\vartheta$ w.r.t. $H$.

Let $\vartheta \in L(H)$ such that $\vartheta^\top \cdot H$ is exponentially special. Then, the Laplace cumulant process of the stochastic integral process $\vartheta^\top \cdot H$ is defined by

$$K(\vartheta) = \kappa(\vartheta) \cdot A,$$

where

$$(4.1) \quad \kappa(\vartheta)_t = \langle \vartheta_t, b_t \rangle + \frac{1}{2}\langle \vartheta_t, c_t \vartheta_t \rangle + \int_{\mathbb{R}^d} (e^{\langle \vartheta_t, x \rangle} - 1 - \langle \vartheta_t, h(x) \rangle) F_t(dx).$$

Analogously to (2.7), it satisfies the following martingale property:

$$(4.2) \qquad \frac{e^{\vartheta^\top \cdot H}}{\mathcal{E}(K(\vartheta))} \in \mathcal{M}_{\mathrm{loc}}(P);$$

cf. Theorems 2.18 and 2.19 in KS. Moreover, $\widetilde{K}$ denotes the logarithmic transform of the cumulant process $K$, that is, $\mathcal{E}(K(\vartheta)) = \exp(\widetilde{K}(\vartheta))$.

THEOREM 4.1. *Let $H$ be an $\mathbb{R}^d$-valued semimartingale with characteristic triplet $\mathbb{T}(H|P) = (B, C, \nu)$. Let $u$ be a vector in $\mathbb{R}^d$. Consider an $\mathbb{R}^d$-valued predictable process $\vartheta$, such that $\vartheta \in L(H)$ and $\vartheta^\top \cdot H$ is exponentially special. Define the measure $P_\vartheta$ via the Radon–Nikodym derivative*

$$\frac{dP_\vartheta}{dP} = \exp(\vartheta^\top \cdot H_T - \widetilde{K}(\vartheta)_T),$$

*assuming that $e^{\vartheta^\top \cdot H - \widetilde{K}(\vartheta)} \in \mathcal{M}(P)$.*

*Then the process $H^u$ with $H^u := u^\top H$, is a 1-dimensional semimartingale with characteristic triplet $\mathbb{T}(H^u|P_\vartheta) = (B^u, C^u, \nu^u)$ of the form*

$$(4.3) \qquad B^u = b^u \cdot A, \qquad C^u = c^u \cdot A, \qquad \nu^u = F^u \otimes A,$$

*where*

$$b^u = u^\top b + u^\top c \cdot \vartheta + \left( h(u^\top x) \frac{e^{\vartheta^\top x}}{1 + W(\vartheta)} - u^\top h(x) \right) * F,$$

$$(4.4) \qquad c^u = u^\top c u,$$

$$F^u(E) = 1_E(u^\top x) \frac{e^{\vartheta^\top x}}{1 + W(\vartheta)} * F, \qquad E \in \mathcal{B}(\mathbb{R} \setminus \{0\}).$$

*Here, $W(\vartheta)_t := \int (e^{\vartheta^\top x} - 1)\nu(\{t\} \times dx)$. (Recall that $h$ denotes a generic truncation function; cf. Remark 2.1.)*

PROOF. We present three proofs of the theorem; the first two proofs reveal interesting relationships between different triplets, while the third proof is "direct" and resembles analogous results for (time-inhomogeneous) Lévy processes; it requires to understand the structure of fixed times of discontinuities.

The structure of the proofs can be represented by the following diagram:

$$(4.5) \quad \mathbb{T}(H|P) \xrightarrow[\substack{(U) \\ (c)}]{\substack{(G) \\ (a)}} \begin{array}{c} \mathbb{T}(H|P_\vartheta) \\ \xrightarrow{(\mathcal{E})} \\ (e) \\ \mathbb{T}(H^u|P) \end{array} \xrightarrow[\substack{(G) \\ (d)}]{\substack{(U) \\ (b)}} \mathbb{T}(H^u|P_\vartheta),$$

where $\xrightarrow{(G)}$ means that we use Girsanov's theorem to calculate the right side triplet from the left side one, $\xrightarrow{(U)}$ means that we use Proposition 2.4 and $\xrightarrow{(\mathcal{E})}$ means that we use the martingale version of the Lévy–Khintchine formula (2.7).

(a) $\mathbb{T}(H|P) \xrightarrow{(G)} \mathbb{T}(H|P_\vartheta)$.

Define the process $Z = (Z_t)_{0 \leq t \leq T}$ via

$$Z := \exp(\vartheta^\top \cdot H - \widetilde{K}(\vartheta)).$$

Clearly, $Z > 0$ a.s., $EZ_T = 1$ and $Z \in \mathcal{M}(P)$ by assumption; cf. KS for conditions. Therefore, the probability measure $P_\vartheta$ defined on $(\Omega, \mathcal{F}, (\mathcal{F}_t)_{0 \leq t \leq T})$ by the Radon–Nikodym derivative $dP_\vartheta = Z_T\, dP$ is equivalent to $P$ ($P_\vartheta \sim P$) and the density process is given by

$$Z_t = \frac{d(P_\vartheta|\mathcal{F}_t)}{d(P|\mathcal{F}_t)} = \exp(\vartheta^\top \cdot H_t - \widetilde{K}(\vartheta)_t), \qquad 0 \leq t \leq T.$$

Moreover, using Theorem 2.19 in KS, we can express $Z$ as follows:

$$(4.6) \qquad Z = \mathcal{E}\left(\vartheta^\top \cdot H^c + \frac{e^{\vartheta^\top x} - 1}{1 + W(\vartheta)} * (\mu - \nu)\right).$$

Now, an application of Girsanov's theorem for semimartingales (JS, Theorem III.3.24) yields that $\mathbb{T}(H|P_\vartheta) = (B^+, C^+, \nu^+)$ where

$$(4.7) \quad \begin{aligned} B^{+i} &= B^i + c^{i\cdot}\beta^+ \cdot A + h^i(x)(Y^+ - 1) * \nu, \\ C^+ &= C, \\ \nu^+ &= Y^+ \cdot \nu, \end{aligned}$$

where, using Lemma 2.20 in KS, we can take the following versions of $\beta^+$ and $Y^+$:

$$\beta^+ = \vartheta \quad \text{and} \quad Y^+ = \frac{e^{\vartheta^\top x}}{1 + W(\vartheta)}.$$

(b) $\mathbb{T}(H|P_\vartheta) \xrightarrow{(U)} \mathbb{T}(H^u|P_\vartheta)$.

Applying Proposition 2.4 to the semimartingale $H$ under the measure $P_\vartheta$ and the vector $u^\top$, we get that $H^u = u^\top H$ is also a $P_\vartheta$-semimartingale with characteristics $\mathbb{T}(H^u|P_\vartheta) = (B^u, C^u, \nu^u)$, where

$$b^u = u^\top b^+ + [h(u^\top x) - u^\top h(x)] * F^+$$
$$= u^\top b + u^\top c \cdot \vartheta + u^\top h(x)(Y^+ - 1) * F$$
(4.8) $$\qquad + [h(u^\top x) - u^\top h(x)]Y^+ * F,$$
$$c^u = u^\top c u,$$
$$F^u(E) = 1_E(u^\top x) * F^+ = 1_E(u^\top x)Y^+ * F, \qquad E \in \mathcal{B}(\mathbb{R} \setminus \{0\}).$$

Therefore, the statement using steps (a) and (b) is proved.

(c) $\mathbb{T}(H|P) \xrightarrow{(U)} \mathbb{T}(H^u|P)$.

Let us denote the triplet $\mathbb{T}(H^u|P) = (B^-, C^-, \nu^-)$; then a direct application of Proposition 2.4 to the process $H$ and the vector $u$ yields that

$$b^- = u^\top b + [h(u^\top x) - u^\top h(x)] * F,$$
(4.9) $$c^- = u^\top c u,$$
$$F^-(E) = 1_E(u^\top x) * F, \qquad E \in \mathcal{B}(\mathbb{R} \setminus \{0\}).$$

(d) $\mathbb{T}(H^u|P) \xrightarrow{(G)} \mathbb{T}(H^u|P_\vartheta)$.

In order to calculate the triplet $\mathbb{T}(H^u|P_\vartheta)$ from the triplet $\mathbb{T}(H^u|P)$, we use Girsanov's theorem for semimartingales (JS, Theorem III.3.24) which states that

$$B^u = B^- + \beta^- c^u \cdot A + h(y)(Y^- - 1) * \nu^-,$$
(4.10) $$C^u = C^-,$$
$$\nu^u = Y^- \cdot \nu^-.$$

Here, $\beta^- = \beta_t^-(\omega)$ and $Y^- = Y^-(\omega; t, y)$ are defined by the following formulas (JS III.3.28):

(4.11) $$\langle Z^c, (H^u)^c \rangle = (Z_- \beta^-) c^u \cdot A$$

and

(4.12) $$Y^- = M_{\mu^{H^u}}^P \left( \frac{Z}{Z_-} \bigg| \widetilde{\mathcal{P}} \right).$$

Note that in (4.12) $\widetilde{\mathcal{P}} = \mathcal{P} \otimes \mathcal{B}(\mathbb{R})$ denotes the $\sigma$-field of predictable sets in $\Omega \times [0, T] \times \mathbb{R}$, and $M^P_{\mu^{H^u}} = \mu^{H^u}(\omega; dt, dy) P(d\omega)$ is the positive measure on $(\Omega \times [0, T] \times \mathbb{R}, \mathcal{F} \otimes \mathcal{B}([0, T]) \otimes \mathcal{B}(\mathbb{R}))$ defined by

$$(4.13) \qquad M^P_{\mu^{H^u}}(W) = E(W * \mu^{H^u})_T$$

for measurable nonnegative functions $W = W(\omega; t, y)$ on $\Omega \times [0, T] \times \mathbb{R}$.

The conditional expectation $M^P_{\mu^{H^u}}(\frac{Z}{Z_-} | \widetilde{\mathcal{P}})$ is, by definition, the $M^P_{\mu^{H^u}}$-a.s. unique $\widetilde{\mathcal{P}}$-measurable function $Y^-$ with the property

$$(4.14) \qquad M^P_{\mu^{H^u}}\left(\frac{Z}{Z_-}U\right) = M^P_{\mu^{H^u}}(Y^- U)$$

for all nonnegative $\widetilde{\mathcal{P}}$-measurable functions $U = U(\omega; t, y)$; cf. JS III.3.16.

We show that in our case, where $Z$ is given by (4.6), we can take the following versions of $\beta^-$ and $Y^-$:

$$(4.15) \qquad \beta^- = \frac{\vartheta^\top u}{|u|^2} \quad \text{and} \quad Y^- = \frac{\exp((\vartheta^\top u/|u|^2)y)}{1 + W(\vartheta)}.$$

Indeed, (4.6) immediately yields that

$$Z^c = Z_- \vartheta^\top \cdot H^c,$$

while we also have that $(H^u)^c = u^\top H^c$. Therefore, using JS I.4.41, III.6.6, we calculate

$$(4.16) \qquad \begin{aligned} \langle Z^c, (H^u)^c \rangle &= \langle Z_- \vartheta^\top \cdot H^c, u^\top H^c \rangle \\ &= Z_- \vartheta^\top \cdot \langle H^c, H^c \rangle \cdot u \\ &= Z_- \vartheta^\top cu \cdot A \\ &= Z_- \vartheta^\top \frac{uu^\top}{|u|^2} cu \cdot A \\ &= \left(Z_- \frac{\vartheta^\top u}{|u|^2}\right) c^u \cdot A, \end{aligned}$$

which yields that we can take $\beta^- = \frac{\vartheta^\top u}{|u|^2}$. Now, for all test functions $U$ and using that $\Delta \widetilde{K}(\vartheta) = \log(1 + W(\vartheta))$, see Theorem 2.18 in KS, we get that

$$\begin{aligned} M^P_{\mu^{H^u}}&\left(\frac{Z}{Z_-}U\right) \\ &= E\left[\int_0^T \int_\mathbb{R} \frac{Z_t(\omega)}{Z_{t-}(\omega)} U(\omega; t, y) \mu^{H^u}(\omega; dt, dy)\right] \\ &= E\left[\sum_{0 \leq t \leq T} e^{\vartheta^\top \Delta H_t(\omega) - \Delta \widetilde{K}(\vartheta)_t} U(\omega; t, u^\top \Delta H_t(\omega)) 1_{\{u^\top \Delta H_t(\omega) \neq 0\}}\right] \end{aligned}$$

$$
\begin{aligned}
(4.17) \\
&= E\bigg[\sum_{0\leq t\leq T} \frac{e^{(\vartheta^\top u/|u|^2)u^\top \Delta H_t(\omega)}}{1+W(\vartheta)_t} U(\omega;t,u^\top \Delta H_t(\omega))1_{\{u^\top \Delta H_t(\omega)\neq 0\}}\bigg] \\
&= E\bigg[\int_0^T \int_\mathbb{R} \frac{e^{(\vartheta^\top u/|u|^2)y}}{1+W(\vartheta)_t} U(\omega;t,y)\mu^{H^u}(\omega;dt,dy)\bigg] \\
&= M^P_{\mu^{H^u}}\bigg(\frac{e^{(\vartheta^\top u/|u|^2)y}}{1+W(\vartheta)_t}U\bigg);
\end{aligned}
$$

therefore, one may take $Y^-$ as in (4.15).

Now, replacing $\beta^-$ and $Y^-$ from (4.15) to (4.10), and using also (4.9), we get that $C^u = C^- = u^\top C u$,

$$
\begin{aligned}
(4.18) \\
1_E(y) * \nu^u &= 1_E(y) * \bigg(\frac{e^{(\vartheta^\top u/|u|^2)y}}{1+W(\vartheta)} \cdot \nu^-\bigg) \\
&= 1_E(u^\top x)\frac{e^{(\vartheta^\top u/|u|^2)u^\top x}}{1+W(\vartheta)} * \nu = 1_E(u^\top x)\frac{e^{\vartheta^\top x}}{1+W(\vartheta)} * \nu
\end{aligned}
$$

and

$$
\begin{aligned}
B^u &= u^\top B + [h(u^\top x) - u^\top h(x)] * \nu + \frac{\vartheta^\top u}{|u|^2}u^\top cu \cdot A \\
(4.19) \quad &+ h(u^\top x)\bigg(\frac{e^{(\vartheta^\top u/|u|^2)u^\top x}}{1+W(\vartheta)_t} - 1\bigg) * \nu \\
&= u^\top B + u^\top c\vartheta \cdot A + \bigg(h(u^\top x)\frac{e^{\vartheta^\top x}}{1+W(\vartheta)_t} - u^\top h(x)\bigg) * \nu.
\end{aligned}
$$

Therefore, the proof using steps (c) and (d) yields the required results.

(e) $\mathbb{T}(H|P) \xrightarrow{(\mathcal{E})} \mathbb{T}(H^u|P_\vartheta)$.

Using the martingale version of the Lévy–Khintchine formula, it is sufficient to show that

$$
(4.20) \quad \frac{e^{zH^u}}{\mathcal{E}(K^u(z))} \in \mathcal{M}_{\text{loc}}(P_\vartheta)
$$

for all $z \in \mathcal{Z}$, where $\mathcal{Z} \subseteq \mathbb{R}$ open, such that

$$
(4.21) \quad 1_{\{|y|>1\}}e^{zy} * \nu^u \in \mathcal{V};
$$

cf. JS III.7.15. Using Proposition III.3.8 in JS, (4.20) is equivalent to

$$
(4.22) \quad Z\frac{e^{zH^u}}{\mathcal{E}(K^u(z))} = \frac{e^{\vartheta^\top \cdot H}}{\mathcal{E}(K(\vartheta))}\frac{e^{zH^u}}{\mathcal{E}(K^u(z))} \in \mathcal{M}_{\text{loc}}(P);
$$

here $K^u$ denotes the Laplace cumulant process associated to the triplet $(B^u, C^u, \nu^u)$. Moreover, using (4.4), condition (4.21) translates to

$$1_{\{|u^\top x|>1\}} e^{(zu+\vartheta)^\top x} * \nu \in \mathcal{V}, \tag{4.23}$$

because $W(\vartheta)$ does not depend on $x$.

Now, the exponent of the numerator in (4.22) is

$$\vartheta^\top \cdot H + zH^u = \vartheta^\top \cdot H + zu^\top H = (\vartheta + zu)^\top \cdot H,$$

which has the *unique* exponential compensator (KS, Lemma 2.15)

$$K(\vartheta + zu)$$

under the measure $P$, for the cumulant defined in (2.5); that is,

$$\frac{e^{(\vartheta+zu)^\top \cdot H}}{\mathcal{E}(K(\vartheta + zu))} \in \mathcal{M}_{\mathrm{loc}}(P). \tag{4.24}$$

Therefore, to complete the proof using step (e), it suffices to show

$$\mathcal{E}(K(\vartheta))\mathcal{E}(K^u(z)) = \mathcal{E}(K(\vartheta + zu)). \tag{4.25}$$

Now, Yor's formula [cf. Jacod (1979), Proposition 6.4] yields that

$$\mathcal{E}(K(\vartheta))\mathcal{E}(K^u(z)) = \mathcal{E}(K(\vartheta) + K^u(z) + [K(\vartheta), K^u(z)]).$$

Hence, (4.25) reduces to showing that

$$K(\vartheta) + K^u(z) + [K(\vartheta), K^u(z)] = K(\vartheta + zu). \tag{4.26}$$

On the right-hand side of (4.26), we have

$$\begin{aligned}K(\vartheta + zu) &= (\vartheta + zu)^\top b \cdot A + \tfrac{1}{2}(\vartheta + zu)^\top c(\vartheta + zu) \cdot A \\ &\quad + (e^{(\vartheta+zu)^\top x} - 1 - (\vartheta + zu)^\top h(x)) * \nu.\end{aligned} \tag{4.27}$$

On the left-hand side of (4.26), we have similarly

$$K(\vartheta) = \vartheta^\top b \cdot A + \tfrac{1}{2}\vartheta^\top c\vartheta \cdot A + (e^{\vartheta^\top x} - 1 - \vartheta^\top h(x)) * \nu; \tag{4.28}$$

the second term on the left-hand side of (4.26), using (4.4), is

$$\begin{aligned}K^u(z) &= zB^u + \frac{1}{2}z^2 C^u + (e^{zy} - 1 - zh(y)) * \nu^u \\ &= zu^\top B + zu^\top c\vartheta \cdot A + z\left(h(u^\top x)\frac{e^{\vartheta^\top x}}{1 + W(\vartheta)} - u^\top h(x)\right) * \nu \\ &\quad + \frac{1}{2}z^2 u^\top cu \cdot A + (e^{zu^\top x} - 1 - zh(u^\top x))\frac{e^{\vartheta^\top x}}{1 + W(\vartheta)} * \nu \\ &= zu^\top B + zu^\top c\vartheta \cdot A + \frac{1}{2}z^2 u^\top cu \cdot A \\ &\quad + \left(\frac{e^{(zu+\vartheta)^\top x} - e^{\vartheta^\top x}}{1 + W(\vartheta)} - zu^\top h(x)\right) * \nu.\end{aligned} \tag{4.29}$$

The last term on the left-hand side of (4.26) is

$$[K(\vartheta), K^u(z)] = \sum_{t \leq \cdot} \Delta K(\vartheta)_t \Delta K^u(z)_t, \tag{4.30}$$

since $K$, $K^u$ are predictable processes of finite variation; cf. JS I.4.53.

Now, we proceed as follows. First, we show that the drift terms on the left- and right-hand side of (4.26) are equal. Then we show that the diffusive terms are equal. Finally, we prove that the jump terms on the left- and right-hand side of (4.26) are equal.

The drift term and the diffusive term are rather easy to handle; indeed from (4.28) and (4.29), we have that the drift term of the LHS of (4.26) is

$$\vartheta^\top b \cdot A + z u^\top b \cdot A = (\vartheta + z u)^\top \cdot A. \tag{4.31}$$

Similarly, the diffusive term of the LHS of (4.26) is

$$\tfrac{1}{2}(\vartheta^\top c \vartheta + z^2 u^\top c u + z u^\top c \vartheta + z \vartheta^\top c u) \cdot A = \tfrac{1}{2}(\vartheta + zu)^\top c(\vartheta + zu) \cdot A, \tag{4.32}$$

since the matrix $C$ is symmetric. Hence, both terms agree with the RHS.

The jump terms are more difficult to manipulate due to the presence of the fixed times of discontinuity for the semimartingale $H$, which entails that the Laplace cumulant process is discontinuous. Regarding the fixed times of discontinuity, using Theorem 2.18 in KS, we have

$$\Delta K(\vartheta)_t = \int (e^{\vartheta^\top x} - 1)\nu(\{t\} \times dx) = W(\vartheta)_t, \tag{4.33}$$

and

$$\Delta K^u(z)_t = \int (e^{zy} - 1)\nu^u(\{t\} \times dy)$$

$$= \int (e^{zu^\top x} - 1)\frac{e^{\vartheta^\top x}}{1 + W(\vartheta)_t}\nu(\{t\} \times dx) \quad \text{(by (4.4))} \tag{4.34}$$

$$= \int (e^{(\vartheta + zu)^\top x} - e^{\vartheta^\top x})\frac{1}{1 + W(\vartheta)_t}\nu(\{t\} \times dx).$$

An important observation here is that $W(\vartheta)$ does not depend on the integrating variable $x$; hence, we can pull it out of the integration, and get

$$\Delta K^u(z)_t = \frac{1}{1 + W(\vartheta)_t} \int (e^{(\vartheta + zu)^\top x} - e^{\vartheta^\top x})\nu(\{t\} \times dx). \tag{4.35}$$

REMARK 4.2. A second important observation is the following: assume for a moment that $\nu$ is a measure of finite variation, that is, $(|x| \wedge 1) * \nu \in \mathcal{V}$, such that the integrals make sense; then

$$
\begin{aligned}
\sum_{t\leq \cdot} \Delta K(\vartheta)_t &= \sum_{t\leq \cdot} \int_{\mathbb{R}^d} (e^{\vartheta^\top x} - 1)\nu(\{t\} \times dx) \\
&= \int_0^\cdot \int_{\mathbb{R}^d} (e^{\vartheta^\top x} - 1)\nu(dt, dx) \\
&= (e^{\vartheta^\top x} - 1) * \nu,
\end{aligned}
\tag{4.36}
$$

since $K$ is a process of finite variation that jumps only at fixed times.

Hence, for the last term in (4.26), we can calculate further using (4.35)

$$
\begin{aligned}
[K(\vartheta), K^u(z)] &= \sum_{t\leq \cdot} \Delta K(\vartheta)_t \Delta K^u(z)_t \\
&= \sum_{t\leq \cdot} W(\vartheta)_t \frac{1}{1+W(\vartheta)_t} \int (e^{(\vartheta+zu)^\top x} - e^{\vartheta^\top x})\nu(\{t\} \times dx) \\
&= \sum_{t\leq \cdot} \frac{(1+W(\vartheta)_t - 1)}{1+W(\vartheta)_t} \int (e^{(\vartheta+zu)^\top x} - e^{\vartheta^\top x})\nu(\{t\} \times dx) \\
&= \sum_{t\leq \cdot} \bigg( \int (e^{(\vartheta+zu)^\top x} - e^{\vartheta^\top x})\nu(\{t\} \times dx) \\
&\quad - \frac{1}{1+W(\vartheta)_t} \int (e^{(\vartheta+zu)^\top x} - e^{\vartheta^\top x})\nu(\{t\} \times dx) \bigg).
\end{aligned}
\tag{4.37}
$$

Now, using the above observations, we can show that the jump terms on the left- and right-hand side of (4.26) are equal. Indeed, we can express the integrals w.r.t. the compensator $\nu$ as sums in (4.28) and (4.29) (see Remark 4.2 for the intuition). Then, we have that the "jump" term on the LHS of (4.26) (denoted by $\mathbb{I}$) is

$$
\begin{aligned}
\mathbb{I} &= \sum_{t\leq \cdot} \int \bigg( e^{\vartheta^\top x} - 1 - \vartheta^\top h(x) + \frac{e^{(zu+\vartheta)^\top x} - e^{\vartheta^\top x}}{1+W(\vartheta)_t} \\
&\quad - zu^\top h(x) + e^{(\vartheta+zu)^\top x} \\
&\quad - e^{\vartheta^\top x} - \frac{1}{1+W(\vartheta)_t}(e^{(\vartheta+zu)^\top x} - e^{\vartheta^\top x}) \bigg)\nu(\{t\} \times dx) \\
&= \sum_{t\leq \cdot} \int (e^{(\vartheta+zu)^\top x} - 1 - (\vartheta+zu)^\top h(x))\nu(\{t\} \times dx) \\
&= (e^{(\vartheta+zu)^\top x} - 1 - (\vartheta+zu)^\top h(x)) * \nu,
\end{aligned}
\tag{4.38}
$$

which equals the corresponding quantity on the RHS of (4.26). This settles the proof using step (e). □

REMARK 4.3. Naturally, we can recover the results about the dual process under the dual measure in dimension one (see Theorem 3.4 in EPS) as a special case of the multidimensional framework. Indeed, assume that $H$ is an $\mathbb{R}$-valued semimartingale such that

$$S = e^H \in \mathcal{M}(P) \quad \text{and} \quad \frac{dP_\vartheta}{dP} = S = e^H,$$

hence $\vartheta = 1$; in this case we have immediately that $K(1) = \widetilde{K}(1) = 0$ and $W(1) = 0$ (cf. Remark 3, page 408 in KS). Moreover, the dual process in the one-dimensional case is $H^u = -H$, hence, $u = -1$.

Then $\mathbb{T}(H'|P') = \mathbb{T}(H^u|P_\vartheta) = (B^u, C^u, \nu^u)$, where

$$B^u = -B - C - h(x)(e^x - 1) * \nu = B',$$
$$C^u = C = C',$$
$$1_E(x) * \nu^u = 1_E(-x)e^x * \nu = 1_E(x) * \nu', \qquad E \in \mathcal{B}(\mathbb{R} \setminus \{0\}),$$

where we have used the antisymmetry property of the truncation function, that is, $h(-x) = -h(x)$.

COROLLARY 4.4. *Let $H$ be a Lévy process, a special and exponentially special semimartingale. Consider $u \in \mathbb{R}^d$ and $\vartheta \in L(H)$ and define the measure $P_\vartheta$ as in Theorem 4.1. Then the process $H^u$ is a 1-dimensional semimartingale with characteristics $\mathbb{T}(H^u|P_\vartheta) = (B^u, C^u, \nu^u)$ of the form*

(4.39)
$$B^u = \int_0^\cdot b_s^u \, ds, \qquad C^u = \int_0^\cdot c_s^u \, ds,$$
$$\nu^u([0, \cdot] \times E) = \int_0^\cdot \int_E F_s^u(dx) \, ds,$$

*where*

(4.40)
$$b_s^u = u^\top b + u^\top c \vartheta_s + \int_{\mathbb{R}^d} u^\top x (e^{\vartheta_s^\top x} - 1) F(dx),$$
$$c_s^u = u^\top c u,$$
$$F_s^u(E) = \int_{\mathbb{R}^d} 1_E(u^\top x) e^{\vartheta_s^\top x} F(dx), \qquad E \in \mathcal{B}(\mathbb{R} \setminus \{0\}).$$

*Note that $H$ is not necessarily a Lévy process under the measure $P_\vartheta$. It remains a Lévy process (PIIS) if $\vartheta$ is deterministic and time-independent; it becomes a time-inhomogeneous Lévy process (PII) if $\vartheta$ is deterministic but time-dependent; it is a general semimartingale if $\vartheta$ is random.*

PROOF. Directly from Theorems 4.1 and II.4.15 and Corollary II.4.19 in JS. □

COROLLARY 4.5. *Let $H$ be an $\mathbb{R}^d$-valued diffusion process that satisfies the stochastic differential equation*

$$(4.41) \qquad dH_t = b(t, H_t)\, dt + \sigma(t, H_t)\, dW_t, \qquad H_0 = 0,$$

*where $b:[0,T] \times \mathbb{R}^d \to \mathbb{R}^d$ and $c := \sigma\sigma^\top$, such that $c:[0,T] \times \mathbb{R}^d \to \mathbb{R}^{d \times d}$ is symmetric and nonnegative definite. The characteristics of $H$ are $(B, C, \nu)$, where obviously $\nu \equiv 0$, and*

$$(4.42) \qquad \begin{aligned} B &= \int_0^\cdot b(s, H_s)\, ds, \\ C &= \int_0^\cdot c(s, H_s)\, ds = \int_0^\cdot \sigma(s, H_s)\sigma(s, H_s)^\top\, ds. \end{aligned}$$

*Consider $u \in \mathbb{R}^d$ and $\vartheta \in \mathbb{R}^d$ (deterministic, for simplicity), and define the measure $P_\vartheta$ as in Theorem 4.1. Then the process $H^u$ is a univariate diffusion process with characteristics $\mathbb{T}(H^u | P_\vartheta) = (B^u, C^u, \nu^u)$ of the form*

$$(4.43) \qquad \begin{aligned} B^u &= u^\top B + u^\top C \vartheta = \int_0^\cdot (u^\top b(s, H_s) + u^\top c(s, H_s)\vartheta)\, ds, \\ C^u &= u^\top C u = \int_0^\cdot (u^\top c(s, H_s) u)\, ds \end{aligned}$$

*and $\nu^u \equiv 0$.*

PROOF. Directly from Theorem 4.1 and JS III.2.19, III.2.23, III.2.27. □

**5. Applications: Models and options.** 1. The payoff of a *swap* option, also coined a "Margrabe" option or option to exchange one asset for another, is

$$(S_T^1 - S_T^2)^+$$

and we denote its value by

$$(5.1) \qquad \mathbb{M}(S^1, S^2) = E[(S_T^1 - S_T^2)^+].$$

The payoff of the *quanto call* and *put* option, respectively, is

$$S_T^1(S_T^2 - K)^+ \quad \text{and} \quad S_T^1(K - S_T^2)^+,$$

and we will use the following notation for the value of the quanto call option

$$(5.2) \qquad \mathbb{QC}(S^1, S^2, K) = E[S_T^1(S_T^2 - K)^+]$$

and similarly for the quanto put option

$$\mathbb{QP}(S^1, K, S^2) = E[S_T^1(K - S_T^2)^+]. \tag{5.3}$$

The different variants of the quanto option traded in Foreign Exchange markets are explained in detail in Musiela and Rutkowski (2005).

The payoff of a *digital* (cash-or-nothing) and a *correlation*, or *quanto*, *digital* option, respectively, is

$$1_{\{S_T > K\}} \quad \text{and} \quad S_T^1 1_{\{S_T^2 > K\}}.$$

Hence, the holder of a correlation digital option receives one unit of the payment asset $(S^1)$ at expiration, if the measurement asset $(S^2)$ ends up in the money. Of course, this is a generalization of the (standard) digital *asset-or-nothing* option, where the holder receives one unit of the asset if it ends up in the money. We denote the value of the digital option by

$$\mathbb{D}(S, K) = E[1_{\{S_T > K\}}] \tag{5.4}$$

and the value of the correlation digital option by

$$\mathbb{CD}(S^1, S^2, K) = E[S_T^1 1_{\{S_T^2 > K\}}]. \tag{5.5}$$

Moreover, we denote the values of the standard call and put options by

$$\mathbb{C}(S, K) = E[(S_T - K)^+] \tag{5.6}$$

and

$$\mathbb{P}(K, S) = E[(K - S_T)^+]. \tag{5.7}$$

THEOREM 5.1. *Assume that the asset price processes evolve as exponential semimartingales according to (3.1)–(3.3), Assumption (ES) is in force and $e^{H^i} \in \mathcal{M}(P)$, $i = 1, 2$. Then we can relate the value of a swap and a plain vanilla option via the following duality:*

$$\mathbb{M}(S^1, S^2) = \mathbb{P}_\vartheta(1, S^u) = \mathbb{C}_\theta(S^v, 1), \tag{5.8}$$

*where the characteristics $(C^u, \nu^u)$ and $(C^v, \nu^v)$ of $H^u = \log S^u$ and $H^v = \log S^v$, respectively, are given by Theorem 4.1 for $\vartheta = (1, 0)^\top$, $u = (-1, 1)^\top$, and $\theta = (0, 1)^\top$, $v = (1, -1)^\top$.*

PROOF. We will use asset $S^1$ as the numéraire asset; if we use asset $S^2$ instead, then we get the duality relationship with a call option. The value of the swap, or "Margrabe," option is

$$\mathbb{M}(S^1, S^2) = E[(S_T^1 - S_T^2)^+]$$

$$= E\left[e^{H_T^1}\left(1 - \frac{S_T^2}{S_T^1}\right)^+\right] = E\left[e^{\langle \vartheta, H_T \rangle}\left(1 - \frac{S_T^2}{S_T^1}\right)^+\right], \tag{5.9}$$

where $\vartheta = (1,0)^\top$. Moreover, $e^{\langle \vartheta, H \rangle} \in \mathcal{M}(P)$ by assumption. Define a new measure $P_\vartheta$ via the Radon–Nikodym derivative

$$\frac{dP_\vartheta}{dP} = e^{\langle \vartheta, H_T \rangle}$$

and the valuation problem (5.9) takes the form

$$\mathbb{M}(S^1, S^2) = E_\vartheta\left[\left(1 - \frac{S_T^2}{S_T^1}\right)^+\right],$$

where we define the process $S^u = (S_t^u)_{0 \leq t \leq T}$ via

(5.10) $$S_t^u = \frac{S_t^2}{S_t^1} = \frac{e^{H_t^2}}{e^{H_t^1}} = e^{\langle u, H_t \rangle} = e^{H_t^u}, \qquad 0 \leq t \leq T,$$

for $u = (-1, 1)^\top$. The triplet of predictable characteristics of the semimartingale $H^u$ is given by Theorem 4.1 for $\vartheta = (1,0)^\top$ and $u = (-1,1)^\top$.

Now, applying Proposition III.3.8 in JS, we obtain that

$$e^{\langle u, H \rangle} \in \mathcal{M}(P_\vartheta) \quad \text{since} \quad e^{\langle u, H \rangle} e^{\langle \vartheta, H \rangle} = e^{H^2} \in \mathcal{M}(P).$$

Therefore, we can conclude that

$$\mathbb{M}(S^1, S^2) = E_\vartheta[(1 - S_T^u)^+]. \qquad \square$$

THEOREM 5.2. *Assume that the asset price processes evolve as exponential semimartingales according to (3.1)–(3.3), Assumption (ES) is in force and $e^{H^i} \in \mathcal{M}(P)$, $i = 1, 2$. Then we can relate the value of a quanto call and a plain vanilla call option via the following duality:*

(5.11) $$\mathbb{QC}(S^1, S^2, K) = \mathbb{C}_\vartheta(S^u, K),$$

*where the characteristics $(C^u, \nu^u)$ of $H^u = \log S^u$ are given by Theorem 4.1 for $\vartheta = (1,0)^\top$ and $u = (0,1)^\top$. An analogous duality result relates the quanto put option and the standard put option.*

PROOF. The value of the quanto call option is

(5.12)
$$\begin{aligned}
\mathbb{QC}(S^1, S^2, K) &= E[S_T^1 (S_T^2 - K)^+] \\
&= E[e^{H_T^1}(S_T^2 - K)^+] \\
&= E[e^{\langle \vartheta, H_T \rangle}(e^{H_T^2} - K)^+] \\
&= E_\vartheta[(e^{H_T^u} - K)^+],
\end{aligned}$$

where $\frac{dP_\vartheta}{dP} = e^{\langle \vartheta, H_T \rangle}$ for $\vartheta = (1,0)^\top$ and $H^u = u^\top H$ for $u = (0,1)^\top$. Hence, the statement is proved. $\square$

REMARK 5.3. Note that

$$e^{H^u} = e^{H^2} \notin \mathcal{M}(P_\vartheta) \quad \text{because} \quad e^{H^u} e^{\langle \vartheta, H \rangle} = e^{H^1 + H^2} \notin \mathcal{M}(P).$$

Hence, this result is a useful computational tool, but cannot serve as a "dual market" theory.

THEOREM 5.4. *Assume that the asset price processes evolve as exponential semimartingales according to (3.1)–(3.3), Assumption (ES) is in force and $e^{H^i} \in \mathcal{M}(P)$, $i = 1, 2$. Then we can relate the value of a correlation digital option and a standard digital option via the following duality:*

(5.13) $$\mathbb{CD}(S^1, S^2, K) = \mathbb{D}_\vartheta(S^u, K),$$

*where the characteristics $(C^u, \nu^u)$ of $H^u = \log S^u$ are given by Theorem 4.1 for $\vartheta = (1, 0)^\top$ and $u = (0, 1)^\top$.*

PROOF. Similar to the proof of Theorem 5.2 and therefore omitted. □

2. In the same framework, we will treat an option that depends on three assets, which will be called, for obvious reasons, a *quanto-swap* option. The payoff of the quanto-swap option is

$$S_T^1 (S_T^2 - S_T^3)^+$$

and can be interpreted as a swap option struck in a foreign currency. Let us denote its value by

$$\mathbb{QS}(S^1, S^2, S^3) = E[S_T^1 (S_T^2 - S_T^3)^+].$$

Here, we will restrict ourselves to semimartingales with *independent increments* (PII).

THEOREM 5.5. *Assume that the asset price processes evolve as exponential semimartingale PIIs according to (3.1)–(3.3), Assumption (ES) is in force and $e^{H^i} \in \mathcal{M}_{\mathrm{loc}}(P)$, $i = 1, 2, 3$. Assume further that $\vartheta^\top H$ is exponentially special for $\vartheta = (1, 1, 0)$ and $\frac{e^{\langle \vartheta, H \rangle}}{\mathcal{E}(K(\vartheta))} \in \mathcal{M}(P)$. Then we can relate the value of a quanto swap option and a standard put (or call) option via the following duality:*

(5.14) $$\mathbb{QS}(S^1, S^2, S^3) = \mathbb{CP}_\vartheta(1, S^u),$$

*where the characteristics $(C^u, \nu^u)$ of $H^u = \log S^u$ are given by Theorem 4.1 for $\vartheta = (1, 1, 0)^\top$ and $u = (0, -1, 1)^\top$. Moreover, the constant $\mathbb{C} := \mathcal{E}(K(\vartheta))$.*

PROOF. Instead of changing measure once using $S^1$ as the numéraire and then once more using either $S^2$ or $S^3$, we will combine $S^1$ and $S^2$ (or $S^3$) directly. We have

$$\begin{aligned}
\mathbb{QS}(S^1, S^2, S^3) &= E[S_T^1(S_T^2 - S_T^3)^+] \\
&= E\left[S_T^1 S_T^2 \left(1 - \frac{S_T^3}{S_T^2}\right)^+\right] \\
&= E[e^{\langle \vartheta, H_T \rangle}(1 - e^{\langle u, H_T \rangle})^+],
\end{aligned} \tag{5.15}$$

where $\vartheta = (1, 1, 0)^\top$ and $u = (0, -1, 1)^\top$.

Clearly, $e^{\langle \vartheta, H_T \rangle} \notin \mathcal{M}(P)$, but since $\vartheta^\top H$ is exponentially special it has a unique exponential compensator; the exponential compensator is given by $K(\vartheta)$ and is deterministic, since $H$ is a PII. Hence, we have

$$\begin{aligned}
\mathbb{QS}(S^1, S^2, S^3) &= \mathcal{E}(K(\vartheta))E\left[\frac{e^{\langle \vartheta, H_T \rangle}}{\mathcal{E}(K(\vartheta))}(1 - e^{\langle u, H_T \rangle})^+\right] \\
&= \mathcal{C}E_\vartheta[(1 - e^{H_T^u})^+],
\end{aligned} \tag{5.16}$$

where $\frac{dP_\vartheta}{dP} = e^{\langle \vartheta, H_T \rangle}$ and $H^u = u^\top H$. Therefore, the triplet $\mathbb{T}(H^u | P_\vartheta)$ can be calculated using Theorem 4.1 with $\vartheta = (1, 1, 0)^\top$ and $u = (0, -1, 1)^\top$, and the statement is proved. □

3. The aim of this final section is to further calculate some explicit examples, especially for processes with jumps. In case the driving semimartingale is continuous, then Theorem 4.1 states that we are dealing again with a continuous semimartingale, usually of the same class, and the characteristics can be calculated very easily. When dealing with semimartingales that exhibit jumps, then one must work a little bit more.

Here, we first revisit the classical result of Margrabe (1978) in the Black–Scholes model, although the same calculations are valid for any continuous semimartingale. The next example involves quanto options and a multidimensional generalization of Merton's jump-diffusion model [cf. Merton (1976)].

The final and more interesting example involves swap and quanto options for multidimensional generalized hyperbolic Lévy processes. The moral of this example can be summarized as follows: when modeling assets by a two-dimensional generalized hyperbolic Lévy model, then the valuation of swap and quanto options is equivalent to the valuation of a call or put option in a one-dimensional generalized hyperbolic Lévy model with suitable parameters.

EXAMPLE 5.6 [Margrabe (1978)]. Consider two assets, $S^1$ and $S^2$, where the dynamics of each asset are

(5.17) $\quad S_t^i = \exp(H_t^i) = \exp(b^i t + \sigma_i W_t^i), \qquad i = 1, 2, 0 \leq t \leq T;$

hence, $H^i$ is a Brownian motion with drift, $i = 1, 2$. In other words, the local or differential characteristics of $H = (H^1, H^2)$ are

$$c = \begin{pmatrix} \sigma_1^2 & \rho\sigma_1\sigma_2 \\ \rho\sigma_1\sigma_2 & \sigma_2^2 \end{pmatrix} \quad \text{and} \quad F \equiv 0,$$

where $\sigma_1, \sigma_2 \geq 0$ and $\rho \in [-1, 1]$ is the correlation coefficient of the Brownian motions $W^1$ and $W^2$, i.e. $\langle W^1, W^2 \rangle = \rho$. Assume, as in Margrabe (1978), that the assets pay no dividends. According to (3.3), the drift characteristic has the form

$$b = -\frac{1}{2}\begin{pmatrix} c_{11} \\ c_{22} \end{pmatrix} = -\frac{1}{2}\begin{pmatrix} \sigma_1^2 \\ \sigma_2^2 \end{pmatrix}.$$

The price of the option to exchange asset $S^1$ for asset $S^2$, according to Theorem 5.1, is equal to the price of a put option with strike 1, on an asset $S^u$ with characteristics $(C^u, \nu^u)$ described by Theorem 4.1 for $\vartheta = (1, 0)^\top$ and $u = (-1, 1)^\top$. Hence, we get that

$$c^u = u^\top c u = \begin{pmatrix} -1 & 1 \end{pmatrix} \begin{pmatrix} \sigma_1^2 & \rho\sigma_1\sigma_2 \\ \rho\sigma_1\sigma_2 & \sigma_2^2 \end{pmatrix} \begin{pmatrix} -1 \\ 1 \end{pmatrix} = \sigma_1^2 + \sigma_2^2 - 2\rho\sigma_1\sigma_2$$

and $F^u \equiv 0$. Therefore, we have recovered the original result of Margrabe [cf. Margrabe (1978), page 179] as a special case in our setting.

Moreover, we have that the drift term $b^u$ of $S^u$ has the form

$$\begin{aligned} b^u &= u^\top b + u^\top c \vartheta \\ &= -\frac{1}{2}\begin{pmatrix} -1 & 1 \end{pmatrix}\begin{pmatrix} \sigma_1^2 \\ \sigma_2^2 \end{pmatrix} + \begin{pmatrix} -1 & 1 \end{pmatrix}\begin{pmatrix} \sigma_1^2 & \rho\sigma_1\sigma_2 \\ \rho\sigma_1\sigma_2 & \sigma_2^2 \end{pmatrix}\begin{pmatrix} 1 \\ 0 \end{pmatrix} \\ &= -\frac{1}{2}(\sigma_1^2 + \sigma_2^2 - 2\rho\sigma_1\sigma_2) = -\frac{1}{2}c^u, \end{aligned}$$

as was expected, since $S^u$ is a $P_\vartheta$-martingale.

EXAMPLE 5.7. Let us consider the following extension of Merton's jump diffusion model: the two assets $S^1$ and $S^2$ are modeled as exponential jump-diffusion processes

(5.18) $\qquad S_t^i = \exp(b_i t + H_t^{i,c} + H_t^{i,d}), \qquad i = 1, 2,$

where (a) the drift terms are determined by the martingale condition (3.3); (b) the continuous martingale parts $H^{i,c}$ are correlated Brownian motions with variance $\sigma_i$ and correlation $\rho\sigma_1\sigma_2$ (as in the previous example); the pure jump parts $H^{i,d}$ are compound Poisson processes with intensity $\lambda_i$, where the jump heights follow, for the sake of brevity, independent normal distributions with variance $\tau_i$ and zero mean. Note that since the Poisson process has finite variation, we can choose the truncation function $h(x) \equiv 0$. Hence, the local characteristics of $H = (H^1, H^2)^\top$ are

$$(5.19) \qquad b = \begin{pmatrix} b_1 \\ b_2 \end{pmatrix}, \qquad c = \begin{pmatrix} \sigma_1^2 & \rho\sigma_1\sigma_2 \\ \rho\sigma_1\sigma_2 & \sigma_2^2 \end{pmatrix},$$

and $F$ has the Lebesgue density $f$ where

$$(5.20) \qquad f(x_1, x_2) = \prod_{i=1}^{2} \frac{1}{\tau_i \sqrt{2\pi}} \exp\left(-\frac{x_i^2}{2\tau_i^2}\right) \cdot \lambda_i.$$

Now, the price of a quanto call option with strike $K$, according to Theorem 5.2, is equal to the price of a call option with the same strike $K$, on an asset $S^u$ with characteristics $(B^u, C^u, \nu^u)$ provided by Theorem 4.1 for $\vartheta = (1,0)^\top$ and $u = (0,1)^\top$. Hence, we can calculate

$$(5.21) \qquad b^u = b_2 + c_{12} = b_2 + \rho\sigma_1\sigma_2,$$

$$(5.22) \qquad c^u = c_{22} = \sigma_2^2;$$

for the Lévy measure, using the independence of the normal variables and completing the square, we have for $y \in \mathbb{R}$, $E \in \mathcal{B}(\mathbb{R}\setminus\{0\})$:

$$1_E(y) * F^u = 1_E(x_2)e^{x_1} * F$$

$$= \int_{\mathbb{R}^2} 1_E(x_2)e^{x_1} \prod_{i=1}^{2}\left(\frac{1}{\tau_i\sqrt{2\pi}} \exp\left(\frac{-x_i^2}{2\tau_i^2}\right)\lambda_i\right) dx_1\, dx_2$$

$$(5.23) \qquad = \lambda_2 \int_E \frac{1}{\tau_2\sqrt{2\pi}} \exp\left(\frac{-x_2^2}{2\tau_2^2}\right) dx_2$$

$$\times \lambda_1 e^{\tau_1^2/2} \underbrace{\int_{\mathbb{R}} \frac{1}{\tau_1\sqrt{2\pi}} \exp\left(-\frac{(x_1 - \tau_1^2)^2}{2\tau_1^2}\right) dx_1}_{=1}$$

$$= \lambda_2 \lambda_1 e^{\tau_1^2/2} \int_E \frac{1}{\tau_2\sqrt{2\pi}} \exp\left(-\frac{x_2^2}{2\tau_2^2}\right) dx_2.$$

Therefore, pricing a quanto option in case the two assets are modeled as jump-diffusions, is equivalent to pricing a call option in a *univariate* jump-diffusion model of the *same class*, with parameters given by (5.21)–(5.23). In particular, jumps occur according to a compound Poisson process with intensity $\lambda^u = \lambda_2 \lambda_1 \exp\frac{\tau_1^2}{2}$, and jump heights are normally distributed with jump variance $\tau_2$ and zero mean.

REMARK 5.8. If we assume that the normal variables in the Lévy measure are correlated, then the statement remains essentially the same. The jump intensity will be different and the normal distribution describing the jumps will have a nonzero mean, but the same variance $\tau_2$.

EXAMPLE 5.9 (GH). Consider two assets that are modeled as (dependent) generalized hyperbolic (henceforth GH) Lévy processes, hence, $H = (H^1, H^2)^\top \sim \mathrm{GH}_2(\lambda, \alpha, \beta, \delta, \mu, \Delta)$; cf. Barndorff-Nielsen (1977). The parameters can take the following values: $\lambda \in \mathbb{R}$, $\alpha, \delta \in \mathbb{R}_{\geq 0}$, $\beta, \mu \in \mathbb{R}^2$, and $\Delta \in \mathbb{R}^{2 \times 2}$ is a symmetric, positive definite matrix; w.l.o.g. we can assume $\det(\Delta) = 1$. Therefore, the triplet of predictable characteristics $(b, c, F)$ is

$$(5.24) \qquad b = \begin{pmatrix} b_1 \\ b_2 \end{pmatrix}, \qquad c \equiv 0,$$

where $b_1, b_2$ are determined by the martingale condition (3.3); the Lévy measure $F$ has the Lebesgue density $f$ where, for $x \in \mathbb{R}^2$,

$$(5.25) \quad f(x) = \frac{e^{\langle \beta, x \rangle}}{\pi \sqrt{\langle x, \Delta^{-1} x \rangle}} \bigg( \int_0^\infty \frac{\sqrt{2y + \alpha^2} K_1(\sqrt{(2y + \alpha^2)\langle x, \Delta^{-1} x \rangle})}{\pi^2 y (J_{|\lambda|}^2(\delta \sqrt{2y}) + Y_{|\lambda|}^2(\delta \sqrt{2y}))} \, dy \\ + \alpha K_1(\alpha \sqrt{\langle x, \Delta^{-1} x \rangle}) 1_{\{\lambda > 0\}} \bigg);$$

cf. Masuda (2004). The limiting case $\delta = 0$, for $\lambda > 0$, corresponds to the bivariate Variance Gamma model and the Lévy measure is

$$(5.26) \qquad f(x) = \frac{\lambda \alpha e^{\langle \beta, x \rangle}}{\pi \sqrt{\langle x, \Delta^{-1} x \rangle}} K_1(\alpha \sqrt{\langle x, \Delta^{-1} x \rangle});$$

cf. Hammerstein (2004). If $\alpha^2 - \langle \beta, \Delta \beta \rangle > 0$, then moments of all orders and the moment generating function exist.

Now, we want to determine the triplet of local characteristics for the univariate Lévy process $H^u$ under the measure $P_\vartheta$, resulting from the duality results for swap and quanto options. According to Theorems 5.1 and 5.2, the characteristics are provided by Theorem 4.1 for $\vartheta = (1, 0)^\top$ and $u = (-1, 1)^\top$ and $\vartheta = (1, 0)^\top$ and $u = (0, 1)^\top$, respectively.

As a first step, we have that the change of probability measure from $P$ to $P_\vartheta$ produces an exponential tilting of the Lévy measure by $\langle \vartheta, x \rangle = x_1$, in both cases. Therefore, under $P_\vartheta$ we have a Lévy process from the same class, with new skewness parameter $\beta_\vartheta = (\beta_1 + 1, \beta_2)$.

The following result determines the parameters of a multidimensional GH distribution under a Radon transformation; the proof follows directly using the characteristic function, see also Lillestøl (2002) for the NIG case. Let $H$ be a random vector such that $H \sim \mathrm{GH}_n(\lambda, \alpha, \beta, \delta, \mu, \Delta)$ and consider the transformation $H^u = u^\top H$, for $u \in \mathbb{R}^n \setminus \{0\}$. Then

$$H^u \sim \mathrm{GH}_1(\lambda^u, \alpha^u, \beta^u, \delta^u, \mu^u),$$

where

$$\lambda^u = \lambda, \qquad \alpha^u = \sqrt{\frac{\alpha^2 - \beta^\top \Delta \beta}{u^\top \Delta u} + \left(\frac{u^\top \Delta \beta}{u^\top \Delta u}\right)^2},$$

$$\beta^u = \frac{u^\top \Delta \beta}{u^\top \Delta u}, \qquad \delta^u = \delta\sqrt{u^\top \Delta u}, \qquad \mu^u = u^\top \mu.$$

Therefore, for the swap option we have that $H^u$ under $P_\vartheta$ follows a univariate GH distribution with parameters

$$\lambda^u = \lambda,$$

$$\alpha^u = \sqrt{\frac{\alpha^2 - (\beta_1 + 1)^2 \delta_{11} - \beta_2^2 \delta_{22} - 2(\beta_1 + 1)\beta_2 \delta_{12}}{\delta_{11} + \delta_{22} - 2\delta_{12}} + (\beta^u)^2},$$

(5.27) $\quad \beta^u = \dfrac{\beta_2 \delta_{22} - (\beta_1 + 1)\delta_{11} - \delta_{12}(\beta_2 - \beta_1 - 1)}{\delta_{11} + \delta_{22} - 2\delta_{12}},$

$$\delta^u = \delta\sqrt{\delta_{11} + \delta_{22} - 2\delta_{12}},$$

$$\mu^u = \mu_2 - \mu_1,$$

hence, the Lévy measure is known. The drift term $b^u$ is determined by the martingale property of $S^u = e^{H^u}$ under the measure $P_\vartheta$.

Similarly, for the quanto option we have that $H^u$ under $P_\vartheta$ follows a univariate GH distribution with parameters

(5.28)
$$\lambda^u = \lambda, \qquad \alpha^u = \sqrt{\frac{\alpha^2 - (\beta_1 + 1)^2(\delta_{11} - \delta_{22}^{-1}\delta_{12}^2)}{\delta_{22}}},$$

$$\beta^u = \beta_2 + \frac{(\beta_1 + 1)\delta_{12}}{\delta_{22}}, \qquad \delta^u = \delta\sqrt{\delta_{22}}, \qquad \mu^u = \mu_2,$$

which also provides the Lévy measure. The drift term is given by

(5.29) $$b^u = b_2 + \int_{\mathbb{R}^2} x_2(e^{x_1} - 1)F(dx).$$

REMARK 5.10. The calculations for quanto options are also valid for correlation digital options.

**Acknowledgment.** We thank an anonymous referee for her/his careful reading of the paper.

# REFERENCES

BARNDORFF-NIELSEN, O. E. (1977). Exponentially decreasing distributions for the logarithm of particle size. *Proc. Roy Soc. London Ser. A* **353** 401–419.

EBERLEIN, E. and PAPAPANTOLEON, A. (2005). Symmetries and pricing of exotic options in Lévy models. In *Exotic Option Pricing and Advanced Lévy Models* (A. KYPRIANON, W. SCHOUTENS AND P. WILMOTT, EDS.) 99–128. Wiley, Chichester. MR2343210

EBERLEIN, E., PAPAPANTOLEON, A. and SHIRYAEV, A. N. (2008). On the duality principle in option pricing: Semimartingale setting. *Finance Stoch.* **12** 265–292. MR2390191

FAJARDO, J. and MORDECKI, E. (2006). Pricing derivatives on two-dimensional Lévy processes. *Int. J. Theor. Appl. Finance* **9** 185–197. MR2215236

GEMAN, H., EL KAROUI, N. and ROCHET, J.-C. (1995). Changes of numéraire, changes of probability measure and option pricing. *J. Appl. Probab.* **32** 443–458. MR1334898

GERBER, H. U. and SHIU, E. S. W. (1996). Martingale approach to pricing perpetual American options on two stocks. *Math. Finance* **6** 303–322. MR1399052

JACOD, J. (1979). *Calcul Stochastique et Problèmes de Martingales. Lecture Notes in Math.* **714**. Springer, Berlin. MR542115

JACOD, J. and SHIRYAEV, A. N. (2003). *Limit Theorems for Stochastic Processes*, 2nd ed. *Grundlehren der Mathematischen Wissenschaften [Fundamental Principles of Mathematical Sciences]* **288**. Springer, Berlin. MR1943877

KALLSEN, J. and SHIRYAEV, A. N. (2001). Time change representation of stochastic integrals. *Teor. Veroyatnost. i Primenen.* **46** 579–585. MR1978671

KALLSEN, J. and SHIRYAEV, A. N. (2002). The cumulant process and Esscher's change of measure. *Finance Stoch.* **6** 397–428. MR1932378

LILLESTØL, J. (2002). Some crude approximation, calibration and estimation procedures for NIG variates. Working paper.

MARGRABE, W. (1978). The value of an option to exchange one asset for another. *J. Finance* **33** 177–186.

MASUDA, H. (2004). On multidimensional Ornstein–Uhlenbeck processes driven by a general Lévy process. *Bernoulli* **10** 97–120. MR2044595

MERTON, R. C. (1976). Option pricing with discontinuous returns. *Bell J. Financ. Econ.* **3** 145–166.

MOLCHANOV, I. and SCHMUTZ, M. (2008). Geometric extension of put-call symmetry in the multiasset setting. Preprint. Available at arXiv: 0806.4506.

MUSIELA, M. and RUTKOWSKI, M. (2005). *Martingale Methods in Financial Modelling*, 2nd ed. *Stochastic Modelling and Applied Probability* **36**. Springer, Berlin. MR2107822

PAPAPANTOLEON, A. (2007). Applications of semimartingales and Lévy processes in finance: Duality and valuation. Ph.D. thesis, Univ. Freiburg.

SATO, K. (1999). *Lévy Processes and Infinitely Divisible Distributions*. Cambridge Univ. Press.

SCHRODER, M. (1999). Changes of numeraire for pricing futures, forwards and options. *Rev. Financ. Stud.* **12** 1143–1163.

V. HAMMERSTEIN, E. A. (2004). Lévy–Khintchine representations of multivariate generalized hyperbolic distributions and some of their limiting cases. Preprint. Univ. Freiburg.

E. EBERLEIN  
DEPARTMENT OF MATHEMATICAL STOCHASTICS  
UNIVERSITY OF FREIBURG  
ECKERSTR. 1  
79104 FREIBURG  
GERMANY  
E-MAIL: eberlein@stochastik.uni-freiburg.de  

A. PAPAPANTOLEON  
FINANCIAL AND ACTUARIAL MATHEMATICS  
VIENNA UNIVERSITY OF TECHNOLOGY  
WIEDNER HAUPTSTRASSE 8  
1040 VIENNA  
AUSTRIA  
E-MAIL: papapan@fam.tuwien.ac.at

A. N. Shiryaev
Steklov Mathematical Institute
Gubkina str. 8
119991 Moscow
Russia
E-mail: albertsh@mi.ras.ru